\documentclass[12pt]{article}

\usepackage[T2A]{fontenc}
\usepackage[cp866]{inputenc}

\usepackage{amsfonts}
\usepackage{amssymb}
\usepackage{amsmath}
\usepackage{amsthm}

\def \le {\leqslant}
\def \ge {\geqslant}

\begin{document}

\begin{Large}
\centerline{\bf Towards BAD conjecture}
\end{Large}
 \vskip+1.5cm \centerline{\bf Moshchevitin N.G. \footnote{ Research is supported by
grants RFFI 06-01-00518 and INTAS 03-51-5070 }} \vskip+1.5cm

\centerline{\bf Abstract.}

For $\alpha, \beta, \delta \in [0,1] ,\,\, \alpha +\beta = 1 $   we consider sets
$$
{\rm BAD}^* (\alpha, \beta ;\delta ) = \left\{\xi = (\xi_1,\xi_2 ) \in [0,1]^2:\,\,\,\inf_{p\in \mathbb{N}} \,\, \max \{ (p\log( p+1))^\alpha
||p\xi_1||,
 (p\log (p+1))^\beta ||p\xi_2||\} \ge \delta \right\}.
$$
We prove that  for different $(\alpha_1,\beta_1), (\alpha_2,\beta_2), \alpha_1 +\beta_1 = \alpha_2 +\beta_2 = 1 $ and $\delta $ small enough
$$
{\rm BAD}^* (\alpha_1, \beta_1 ;\delta ) \bigcap {\rm BAD}^* (\alpha_2, \beta_2 ;\delta ) \neq \varnothing .
$$
 Our result is based on A. Khintchine's construction and an original method  due to Y. Peres and W. Schlag.
 \vskip+1.5cm

{\bf 1.\,\, Introduction.}

For $\alpha, \beta \in [0,1] $ under condition $ \alpha +\beta = 1 $ and $ \delta \in (0,1/2)$
 we consider the sets
$$
{\rm BAD}(\alpha, \beta ;\delta ) = \left\{\xi = (\xi_1,\xi_2 ) \in [0,1]^2:\,\,\,\inf_{p\in \mathbb{N}} \,\, \max \{ p^\alpha ||p\xi_1||,
 p^\beta ||p\xi_2||\} \ge \delta \right\}
$$
(here $||\cdot ||$ denotes the distance to the nearest integer) and
$$
{\rm BAD}(\alpha, \beta  ) = \bigcup_{\delta > 0} {\rm BAD}(\alpha, \beta ;\delta )
 .
$$
  In \cite{SCH} W.M. Schmidt conjectured that for any $\alpha_1,\alpha_2, \beta_1,\beta_2 \in [0,1],\, \alpha_1 +\beta_1 =\alpha_2+\beta_2 =1$
  the intersection
  $$
{\rm BAD}(\alpha_1, \beta_1  ) \bigcap {\rm BAD}(\alpha_2, \beta_2  )
$$
is not empty. This conjecture is still open. Up to now it is only known that for any  $\alpha, \beta \in [0,1],\,\, \alpha +\beta = 1 $ the set
$$
{\rm BAD}(1,0  ) \bigcap {\rm BAD}(\alpha , \beta  ) \bigcap {\rm BAD}(0,1  )
$$
is not empty. In fact  A. Pollington and S. Velani \cite{PV} proved that this intersection has Hausdorff dimension equal to 2. Some
multidimensional generalizations are due to S. Kristensen, R. Thorn, S. Velani \cite{KTV} and R. Akhunzhanov \cite{RA}. We would like to note
that the papers \cite{PV},\cite{KTV},\cite{RA} mentioned above use H. Davenport's approach \cite{D1},\cite{D2}.

 In the present paper we consider sets
$$
{\rm BAD}^* (\alpha, \beta ;\delta ) = \left\{\xi = (\xi_1,\xi_2 ) \in [0,1]^2:\,\,\,\inf_{p\in \mathbb{N}} \,\, \max \{ (p\log( p+1))^\alpha
||p\xi_1||,
 (p\log (p+1))^\beta ||p\xi_2||\} \ge \delta \right\}.
$$
Obviously
$$
{\rm BAD} (\alpha, \beta ;\delta ) \subseteq {\rm BAD}^* (\alpha, \beta ;\delta ).
$$
As  the series $ \sum_{p=1}^\infty \frac{1}{p\log (p+1)}$ diverges one can easy see that the set
$$
{\rm BAD}^*(\alpha, \beta  ) = \bigcup_{\delta > 0} {\rm BAD}^*(\alpha, \beta ;\delta )
$$
has Lebesgue  measure equal to zero.

The aim of this paper is to prove that the intersection of the sets ${\rm BAD}^*(\alpha, \beta ;\delta )$ is not empty.

 {\bf Theorem 1.}\,\, {\it
Let $\alpha_1,\alpha_2, \beta_1,\beta_2 \in [0,1],\, \alpha_1 +\beta_1 =\alpha_2+\beta_2 =1$ and $ 0<\delta \le 2^{-20} .$
 Then the set
$$
{\rm BAD}(0, 1 ;\delta )\bigcap {\rm BAD}^*(\alpha_1, \beta_1 ;\delta ) \bigcap {\rm BAD}^*(\alpha_2, \beta_2 ;\delta )
$$
is not empty.}

The proof is based on old construction due to A. Khintchine \cite{H1} and recent arguments due to Y. Peres and W. Schlag \cite{P}.
 It is possible to prove by our method that the intersection of any finite collection of the sets of the type $
{\cal B}_j =
 \bigcup_{\delta >0}{\rm BAD}^*(\alpha_j, \beta_j ;\delta )$
 with $\alpha_j +\beta_j = 1, j = 1,...,d$ has positive Hausdorff dimension. Also our two-dimensional result admit multidimensional
 generalizations.

{\bf 2.  Bohr sets.}

 For real  $\xi $ define
$$
H_\xi^{(\beta)} (p,q) = \left\{ x\in \mathbb{N}:\,\,\, p< x\le
q,\,\, ||x \xi ||\le \frac{\delta}{(p\log (p+1))^\beta }\right\}
$$
and
$$
K_\xi^{(\beta)} (p,q) = \left\{ x\in \mathbb{N}:\,\,\, p< x\le
q,\,\, ||x \xi ||\le \frac{\delta}{(x\log (x+1))^\beta }\right\}
\subseteq H_\xi^{(\beta)} (p,q).
$$

We define $\xi \in \mathbb{R}$ to be is a $\delta$-badly approximable number if $ \inf_{p\in \mathbb{N}} \,\,  p ||p\xi|| \ge \delta$. For
$\delta$-badly approximable number $\xi$ and any real number $\eta$ the pair $(\eta, \xi )$ belongs to $ {\rm BAD}(0, 1 ;\delta )$.

{\bf  Lemma 2. } \,\, {\it For a $\delta$-badly approximable number $\xi$ the following inequality is valid:
$$
\# H_\xi^{(\beta)}(p,2p)\le \frac{(24\delta + 2)p^\alpha}{(\log
(p+1) )^\beta }.
$$}

Proof. The number of elements in $H_\xi^{(\beta)} (p,2p)$ is
bounded by the number of integer points in the region
$$
\Omega (p) =
 \left\{ (x,y ) \in \mathbb{R}^2:\,\,\, 0\le x\le 2p,\,\, |x \xi  -y|\le \frac{\delta}{(p\log (p+1))^\beta }\right\}.
$$
We should note that $(0,0) \in \Omega (p)$ but $(0,0) \not \in
H_\xi^{(\beta)}(p,2p)$. Hence the number of integer points in
$H_\xi^{(\beta)}(p,2p)$ does not exceed the number of integer
points in $ \Omega (p)$ minus one.

 For  The measure of  $\Omega (p)$ we have
\begin{equation}
\mu \left(\Omega (p)\right) = \frac{4\delta p^\alpha}{(\log
(p+1))^\beta }  . \label{mu}
\end{equation}

Let there exists an integer primitive point  $(p_0,a_0)\in \Omega
(p), p_0 \ge 1$ (otherwise $ \# H_\xi^{(\beta)} (p,2p)=0$). As $
\xi $ is a $\delta$-badly approximable number we have
$$
\frac{\delta}{p_0}\le |p_0\xi  - a_0|\le \frac{\delta}{(p\log (p+1))^\beta }
$$
and $ p_0 \ge (p\log (p+1))^\beta$. Now we consider two cases.

In the {\bf first case} we suppose that all integer points in $\Omega (p)$ are of the form $ \lambda (p_0,a_0)$  with integer $\lambda$.
 Then the number of such points does not exceed
 $[2p/p_0]+1\le 2p/ (p\log (p+1))^\beta+ 1= \frac{2p^\alpha}{(\log (p+1) )^\beta }+1$ and lemma follows.

In the {\bf second case} the convex hull of all integer points in the region $\Omega (p)$ is a convex polygon $\Pi$ with integer vertices. Let
$m$ be the number of integer points in $\Omega (p)$. Then $ m = m_1 +m_2\ge 3$ where $m_1$ is the number of integer points inside $\Pi$ and $
m_2$ is the number of integer points on the boundary of polygon $\Pi$.
 According to Pick's formula
we have
$$
\frac{m}{6} \le m_1+ \frac{m_2}{2} - 1 = \mu (\Pi ) \le \mu (\Omega ).
$$
 Now lemma follows from (\ref{mu}).

{\bf  Corollary 3. } \,\, {\it Put $q = \left[\frac{p^2}{\delta}
\log  \frac{p^2}{\delta}\right]+1$ and let $\delta < 1/24$. Then
for   a $\delta$-badly approximable number $\xi$ and $\alpha \in
[0,1]$ one has
$$
\sum_{x \in K_\xi^{(\beta)} (p,q )} \frac{1}{(x\log (x+1))^\alpha}
\le 2^6(1+\log (1/\delta )).
$$}

Proof.
$$
\sum_{x \in K_\xi^{(\beta)} (p,q)} \frac{1}{(x\log
(x+1))^\alpha}\le \sum_{\nu = 0}^{[\log (q/p))]} \sum_{x \in
H_\xi^{(\beta)} (2^\nu p, 2^{\nu+1}p)}\frac{1}{(x\log
(x+1))^\alpha} \le \sum_{\nu = 0}^{[\log (q/p)]} \frac{ \# H_\beta
(2^{\nu}p,2^{\nu+1}p)}{(2^\nu p\log (p+1))^\alpha }.
$$
Applying  Lemma 2 we obtain
$$
\sum_{x \in K_\xi^{(\beta)} (p,q)} \frac{1}{(x\log
(x+1))^\alpha}\le \frac{4(24\delta+2)}{\log (p+1)} \sum_{\nu =
0}^{[\log (q/p)]} 1 \le 2^{6}(1+\log (1/\delta ))  ,
$$
and Corollary 2 follows.

{\bf  3. Sets of reals.}

For integers $ 1\le x, 0\le y\le x$  define
\begin{equation}
E_\alpha (x,y) = \left( \frac{y}{x} -\frac{\delta}{x^{1+\alpha}(\log (x+1))^\alpha}, \frac{y}{x} +\frac{\delta}{x^{1+\alpha}(\log (x+1))^\alpha}
\right),
 \,\,\,
E_\alpha (x) =\bigcup_{y =0}^{x} E_\alpha (x,y)\bigcap [0,1]. \label{E}
\end{equation}
 Define
\begin{equation}
l(0,\alpha ) = 0,\,\,\,
 l(x,\alpha ) = \left[\frac{\log (x^{1+\alpha }(\log (x+1))^\alpha/2\delta )}{\log
2}\right] ,\,\, x \in \mathbb{N} . \label{L}
\end{equation}
 Each
segment form the union $E_\alpha (x)$ from (\ref{E}) can be covered by a dyadic interval of the form
$$
\left( \frac{b}{2^{l(x,\alpha ) }}, \frac{b+z}{2^{l(x,\alpha ) }}\right),\,\,\, z = 1,2 .$$

Let $A_\alpha (x)$ be the smallest union of all such dyadic segments which cover the whole set  $E_\alpha (x)$. Define
$$
A^c_\alpha (x) = [0,1] \setminus A_\alpha (x). $$ Then
$$
A^c_\alpha (x)
 = \bigcup_{\nu = 1}^{\tau_\nu } I_\nu
$$
where closed  segments $I_\nu $ are of the form
 \begin{equation}
\left[ \frac{a}{2^{l(x,\alpha )}}, \frac{a+1}{2^{l(x,\alpha )}}\right] ,\,\,\, a\in \mathbb{Z}. \label{I}
\end{equation}

Now we take  to be a $\delta$-badly approximable number  and  reals $ \alpha_1, \alpha_2, \beta_2, \beta_2$ from our Theorem 1. Let for
convenience
$$
\max \{ \alpha_1, \alpha_2 \} = \alpha_1.
$$
Define $B_0 = [0,1]$. For $ q \ge 1$ we consider  sets
$$
B_q = \left(\bigcap_{x \in K_\xi^{(\beta_1)} (0,q) } A^c_{\alpha_1} (x) \right) \bigcap \left(\bigcap_{x \in K_\xi^{(\beta_2)} (0,q) }
A^c_{\alpha_2} (x) \right).
$$
 For $ q \ge 0 $ these sets can be represented in the form
$$
B_q = \bigcup_{\nu = 1}^{T_q } J_\nu \bigcup B_q'
$$
where segments $J_\nu $ are of the form
 \begin{equation}
\left[ \frac{b}{2^{l(q,\alpha_1 )}}, \frac{b+1}{2^{l(q,\alpha_1 )}}\right] ,\,\,\, b\in \mathbb{Z}  \label{J}
\end{equation}
and $B_q'$ consist of  only finitely many points.

We should note that for any $q$ the set $B_q$ is a closed set and the sequence of these sets is nested:
$$
B_1 \supseteq B_2\supseteq \cdots \supseteq B_q \supseteq \cdots .
$$
Hence to prove Theorem 1 it is sufficient to find an infinite sequence of values of $q$ for which the sets $B_q$ are nonempty.

{\bf  3. Lower bounds for $\mu (B_q)$.}

 {\bf  Lemma 4. }
\,\, {\it Let $ \mu (B_q )   \neq 0$. Then for every $i \in \{ 1,2\}$ and for
\begin{equation}
x \ge \frac{q^2}{\delta} \log  \frac{q^2}{\delta} \label{LOG}
\end{equation}
 one has
 $$
\mu \left( B_q  \bigcap A_{\alpha_i} (x) \right) \le   \frac{8 \delta \mu \left( B_q\right)}{(x\log(x+1))^{\alpha_i}}.
$$}

Proof. Write
$$
B_q = \bigcup_{\nu = 1}^{T_q } J_\nu  \bigcup B_q' ,\,\,\, A_{\alpha_i} (x) =\bigcup_{k=1}^{\tau_x} I_k,
$$
where $J_\nu $ are of the form (\ref{J}) and $I_i$ are of the form (\ref{I}) with $\alpha = \alpha_i$.
  Let $I_k \cap J_\nu   \neq \varnothing$.
Remember that $E_{\alpha_i}(x) \subseteq A_{\alpha_i}(x)$
   Then for some natural $y$ we have
$$
\frac{y}{x} \in \left[ \frac{b}{2^{l(q,\alpha_i )}} - \frac{1}{2^{l(x,\alpha_i)}}, \frac{b+1}{2^{l(q,\alpha_i
)}}+\frac{1}{2^{l(x,\alpha_i)}}\right]
  $$
 The quantity $y$ here can take not more than
 $$
 W = \left[   \left( \frac{1}{2^{l(q,\alpha_i )} } + \frac{2}{2^{l(x,\alpha_i)}} \right) x\right]
+ 1
$$
values. Now
$$
\mu \left( J_\nu  \bigcap  A_{\alpha_i}(x)\right) \le \mu\left( I_k\right) W
$$
and
$$
\mu \left( B_q  \bigcap  A_{\alpha_i}(x)\right) \le \mu\left( I_k\right) W T_q = \mu \left(B_q \right) \times \frac{2^{l(q,\alpha_i
)}}{2^{l(x,\alpha_i )}} \times \left( \frac{3x}{2^{l(q,\alpha_i )}}+1\right)= \mu \left(B_q \right) \times  \left( \frac{3x}{2^{l(x,\alpha_i
)}}+\frac{2^{l(q,\alpha_i )}}{2^{l(x,\alpha_i )}}\right)
 .
$$
From  (\ref{L}) we see that
$$
\frac{x}{2^{l(x,\alpha_i )}} \le \frac{  2\delta  }{(x\log (x+1))^{\alpha_i}}.
$$
From (\ref{LOG}) it follows that $ x/\log (x+1) \ge  {q^2}/{\delta} $. Hence
$$
\frac{2^{l(q,\alpha_i )}}{2^{l(x,\alpha_i )}} \le 2\times \frac{q^{1+\alpha_i }}{x^{1+\alpha_i }}\le \frac{2}{(x\log (x+1))^{\alpha_i}} \times
\frac{q^2}{x/\log(x+1)}\le \frac{2\delta}{(x\log (x+1))^{\alpha_i}}.
$$
Now Lemma 4 follows.

{\bf  Lemma 5. } \,\, {\it Let  $0<\delta \le 2^{-20}$. Let
 $q_1 =\left[\frac{q^2}{\delta} \log  \frac{q^2}{\delta}  \right]+1,\,\, q_2 = \left[\frac{q_1^2}{\delta} \log
\frac{q_1^2}{\delta} \right]+1$. Let we know that
\begin{equation}
\mu \left( B_{q_1}\right) \ge \frac{\mu \left( B_{q}\right)}{2}
>0.
\label{supp} \end{equation}
 Then
$$
\mu \left( B_{q_2}\right) \ge \frac{\mu \left( B_{q_1}\right)}{2}.
$$}

Proof. As
$$
B_{q_2} = B_{q_1} \setminus\left(
 \left(\bigcup_{x \in K_\xi^{(\beta_1)} (q_1,q_2) } A_{\alpha_1} (x) \right) \bigcup
 \left(\bigcup_{x \in K_\xi^{(\beta_2)} (q_1,q_2) }
A_{\alpha_2} (x) \right)\right)
$$
we see that
\begin{equation}
\mu \left( B_{q_2}\right) \ge \mu \left( B_{q_1}\right) - \sum_{i=1}^2 \sum_{x \in K_\xi^{(\beta_i)} (q_1,q_2) } \mu \left( B_{q_1}\cap
A_{\alpha_i} (x) \right). \label{A} \end{equation}

 But
$$
B_{q_1}\cap A_{\alpha_i} (x)\subseteq B_{q}\cap A_{\alpha_i} (x) .$$ Applying Lemma 4 we have
$$
\mu \left(  B_{q_1}\cap A_{\alpha_i} (x) \right)\le \mu \left(  B_{q}\cap A_{\alpha_i} (x) \right)\le \frac{8 \delta \mu \left(
B_q\right)}{(x\log(x+1))^{\alpha_i}}.
$$
To continue this estimate we use (\ref{supp}) and obtain
\begin{equation}
\mu \left(  B_{q_1}\cap A_{\alpha_i} (x) \right)\le   \frac{16 \delta \mu \left( B_{q_1}\right)}{(x\log(x+1))^{\alpha_i}}. \label{B}
\end{equation}
Then (\ref{A},\ref{B})   lead to the inequality
$$
\mu \left( B_{q_2}\right) \ge \mu \left( B_{q_1}\right) \left( 1
 -
16 \delta  \times
 \sum_{i=1}^2 \sum_{x \in K_\xi^{(\beta_i)} (q_1,q_2) }\frac{1}{(x\log(x+1))^{\alpha_i}} \right) .
$$
For the last two inner sums we use Corollary 3 with $p=q_1, q=q_2$. We should take into account that from
 the condition $\delta < 2^{-20}$
we have $2^{9}\delta (1+\log (1/\delta ) )\le 1/2$. Now Lemma 5 follows.

 Theorem 1 follows from Lemma 5.
Define $ q_0 = 0, q_{\nu+1}  =\left[\frac{q_\nu^2}{\delta} \log  \frac{q_\nu^2}{\delta}  \right]+1$.
 As $\mu (B_0) = 1$ from Lemma 5 we have
 $\mu (B_{q_\nu}) \ge 2^{-\nu}$. It means that $\bigcap_{q\in \mathbb{N}} B_q \neq \varnothing$.

\end{document}